# INNER PRODUCTS AND BANACH ALGEBRA STRUCTURES ON BICOMPLEX NUMBERS AND THEIR ASSOCIATED SPACES

PRABHAT KUMAR*, FAHED ZULFEQARR, AMIT UJLAYAN, ANJALI, AKHIL PRAKASH, AND M. YOUNUS BHAT

ABSTRACT. In this paper, we introduce various types of inner products and norms on the bicomplex number system $\mathbb{C}_2$, the bicomplex vector space ${\mathbb{C}_2}^n$, the space of bicomplex matrices ${C_2}^{m\times n}$, and the space of bicomplex polynomials $\mathbb{C}_2[\xi]_n$. We investigate the relationships among these inner products and norms, and establish several results. Furthermore, we prove that $\mathbb{C}_2$ and ${\mathbb{C}_2}^n$ are Banach algebras and Hilbert spaces. These results provide a unified framework for the study of inner product structures and normed linear spaces over bicomplex numbers and their associated spaces.

## INTRODUCTION

The development of new number systems has long been a cornerstone of mathematical progress, providing deeper frameworks for analysing and modelling complex phenomena. Among these, the bicomplex number system has recently garnered significant attention due to its ability to generalise and extend classical complex analysis. The bicomplex number system extends the complex number system by introducing a second imaginary unit (see [3, 10, 11]), offering a richer algebraic and analytic structure compared to conventional complex numbers.

In bicomplex analysis, functions depend simultaneously on two complex variables, allowing exploration of higher-dimensional behaviours within a unified algebraic setting (see [7, 10, 12]). This dual-complex structure reveals subtle interactions and properties that do not arise in standard complex function theory, thus opening new avenues for both theoretical investigations and practical applications.

Notably, bicomplex numbers are closely related to quaternions and Clifford algebras, situating them within the broader family of hypercomplex systems (see [6, 13]). This connection has motivated research across diverse domains, including wave propagation, electromagnetic theory, quantum mechanics, and signal processing, where multidimensional and anisotropic effects are naturally modelled within the bicomplex framework. Over the years, numerous researchers [1, 2, 3, 4, 8, 9] have made significant contributions to bicomplex analysis and algebra. Their studies have established several important algebraic and analytical results, laying a solid foundation for further investigations in this area.

In this context, the study of inner product structures on bicomplex spaces $\mathbb{C}_2$ and their associated normed and Hilbert-like spaces has become increasingly important. Defining and analysing various inner products compatible with the bicomplex algebra not only enriches the theoretical landscape but also enhances the potential for applications in functional analysis, operator theory, and mathematical physics.

Motivated by these considerations, this paper introduces and examines several distinct inner products on $\mathbb{C}_2$ and its related vector spaces. We present several norm structures, investigate their properties in relation to the conventional norm, and discuss implications for the construction of normed algebras and possible extensions to Hilbert spaces. By doing so, we seek to highlight the structural richness and analytical potential of bicomplex space, and to motivate further studies in this promising direction.

The present paper is organised into two sections. The first recalls some necessary preliminaries and notations, while the second is devoted to establishing our main results.

## 1. Basics and Fundamentals

In this section, we recall some fundamental concepts that will be employed throughout the paper. The bicomplex-theoretic notations and conventions of [3, 10] are adopted consistently, often without explicit mention. The bicomplex number system was first introduced by Corrado Segre [14] in 1892, who demonstrated that it forms a commutative algebra over the real numbers. The bicomplex numbers are a natural extension of complex numbers, and hyperbolic numbers are a special case of bicomplex numbers when $x_1 = 0 = x_2$, in (1). The algebraic properties of bicomplex numbers are also important and useful (see [11]).

The algebra of bicomplex numbers is commutative but includes zero divisors, which implies that not all nonzero elements possess multiplicative inverses. On the other hand, the algebra of quaternions is a division algebra, meaning every nonzero element is invertible, although the multiplication operation is non-commutative. The commutativity gives major differences between the algebras of quaternions and bicomplex numbers.

**Definition 1.1** (**Bicomplex number** [10]:). A number, which is expressible in the form of

$$\xi = x_1 + i_1 x_2 + i_2 x_3 + i_1 i_2 x_4, \tag{1}$$

is called a bicomplex number, where $x_1, x_2, x_3$, and $x_4$ are real numbers, ${i_1}^2 = {i_2}^2 = -1$, and $i_1 i_2 = i_2 i_1$. Furthermore, every bicomplex number $\xi$ can be uniquely expressed in the following three forms: the real form, the complex form, and the idempotent form, as follows:

$$\begin{aligned} \xi &= x_1 + i_1 x_2 + i_2 x_3 + i_1 i_2 x_4 \rightarrow \text{ Real form of } \xi, \\ \xi &= z_1 + i_2 z_2 \rightarrow \text{ Complex form of } \xi, \\ \xi &= \xi^- e_1 + \xi^+ e_2 \rightarrow \text{ Idempotent form of } \xi, \end{aligned}$$

where $z_1 = x_1 + i_1 x_2$, $z_2 = x_3 + i_1 x_4$, and $\xi^-$ and $\xi^+$ are the idempotent components of $\xi$ (see [3, 15]). The notations $\mathbb{C}_2$, $\mathbb{C}_1$, and $\mathbb{C}_0$ are used for the set of all bicomplex numbers, the set of all complex numbers and the set of all real numbers, respectively.

Zero divisors are found in $\mathbb{C}_2$, therefore it does not form a field structure, but $\mathbb{C}_2$ forms an algebra structure over $\mathbb{C}_1$. $\mathbb{C}_2$ contains exactly four idempotent elements namely $0, 1, e_1$, and $e_2$, where $e_1$ and $e_2$ are the nontrivial idempotent elements and

$$e_1 = \frac{(1 + i_1 i_2)}{2} \text{ and } e_2 = \frac{(1 - i_1 i_2)}{2}.$$

Moreover, $e_1 e_2 = 0 = e_2 e_1$, $(e_i)^n = e_i$, for all $i = 1, 2$, $e_1 + e_2 = 1$, and $e_1 - e_2 = i_1 i_2$.

The author, in [14], has shown that a bicomplex number $\xi = z_1 + i_2 z_2$ can be written uniquely as the complex combination of the elements $e_1$ and $e_2$ in the following way:

$$\xi = (z_1 - i_1 z_2) e_1 + (z_1 + i_1 z_2) e_2, \tag{2}$$

where complex combinations $(z_1 - i_1 z_2)$ and $(z_1 + i_1 z_2)$, in (2), are called idempotent components of $\xi$ and denoted by $\xi^-$ and $\xi^+$, respectively. Thus the bicomplex number $\xi = z_1 + i_2 z_2$ can also be written as $\xi = \xi^- e_1 + \xi^+ e_2$, where $\xi^- = z_1 - i_1 z_2$ and $\xi^+ = z_1 + i_1 z_2$ (see [3, 15]).

**Remark 1.2.** [10] The structure $(\mathbb{C}_2, +, \cdot, \times, )$ is an algebra over $\mathbb{C}_1$ and $\mathbb{C}_0$. The operations addition $(+)$, scalar multiplication $(\cdot)$, and multiplication $(\times)$ are perfectly compatible with the idempotent components of bicomplex numbers and the scalars i.e., if $\xi = z_1 + i_2 z_2, \eta = w_1 + i_2 w_2$ and $\alpha \in \mathbb{C}_1$, then

$$\begin{aligned} \xi + \eta &= \left(\xi^- + \eta^-\right) e_1 + \left(\xi^+ + \eta^+\right) e_2, \\ \xi \times \eta &= \left(\xi^- \eta^-\right) e_1 + \left(\xi^+ \eta^+\right) e_2, \\ \alpha \cdot \xi &= \left(\alpha \xi^-\right) e_1 + \left(\alpha \xi^+\right) e_2. \end{aligned}$$

Furthermore, the multiplication of two bicomplex numbers $\xi$ and $\eta$ in complex form is given by

$$\xi \times \eta = (z_1 + i_2 z_2) \times (w_1 + i_2 w_2) = (z_1 w_1 - z_2 w_2) + i_2 (z_1 w_2 + z_2 w_1).$$

**Definition 1.3** (**Standard norm on** $\mathbb{C}_2$ [10]:). Let $\xi$ be an arbitrary element in $\mathbb{C}_2$ such that

$$\xi = x_1 + i_1x_2 + i_2x_3 + i_1i_2x_4 = z_1 + i_2z_2.$$

Then, the norm of $\xi$ is defined in the following equivalent forms:

$$\|\xi\| = \left(|x_1|^2 + |x_2|^2 + |x_3|^2 + \left|x_4\right|^2\right)^{\frac{1}{2}} \rightarrow \text{Real form of } \|\xi\|, \tag{3}$$

$$\|\xi\| = \left(|z_1|^2 + |z_2|^2\right)^{\frac{1}{2}} \rightarrow \text{Complex form of } \|\xi\|, \tag{4}$$

$$\|\xi\| = \left(\frac{\left|\xi^-\right|^2 + \left|\xi^+\right|^2}{2}\right)^{\frac{1}{2}} \rightarrow \text{Idempotent form of } \|\xi\|. \tag{5}$$

This norm $\|\bullet\|$, on $\mathbb{C}_2$, is called the standard norm and these equivalent expressions (3), (4), and (5) are referred to as the real form, complex form, and idempotent form of the norm of $\xi$, respectively. Moreover

$$\|\xi \times \eta\| \leq \sqrt{2}\, \|\xi\| \cdot \|\eta\|. \tag{6}$$

**Remark 1.4.** [10] Due to the inequality (6), $\mathbb{C}_2$ forms a modified Banach algebra under the standard norm.

**Definition 1.5** (**Conjugate of a bicomplex number**[5]:). Analogous to the concept of the conjugate of a number in $\mathbb{C}_1$, there are three types of conjugates of a bicomplex number $\xi = x_1 + i_1x_2 + i_2x_3 + i_1i_2x_4 = z_1 + i_2z_2$. The $i_1$-conjugate, $i_2$-conjugate, and $i_1i_2$-conjugate are denoted by $\overline{\xi}$, $\widetilde{\xi}$, and $\xi^\sharp$, respectively, and defined as follows:

$$\begin{aligned}
\overline{\xi} &= (x_1 - i_1x_2) + i_2(x_3 - i_1x_4),\\
\widetilde{\xi} &= (x_1 + i_1x_2) - i_2(x_3 + i_1x_4),\\
\xi^\sharp &= (x_1 - i_1x_2) - i_2(x_3 - i_1x_4).
\end{aligned}$$

It is very easy to see that $\overline{\xi} = (\overline{z_1} + i_2\overline{z_2}) = (\overline{\xi^+})e_1 + (\overline{\xi^-})e_2$, $\widetilde{\xi} = (z_1 - i_2z_2) = (\xi^+)e_1 + (\xi^-)e_2$, and $\xi^\sharp = (\overline{z_1} - i_2\overline{z_2}) = (\overline{\xi^-})e_1 + (\overline{\xi^+})e_2$.

**Remark 1.6.** [1] If $\xi = x_1 + i_1x_2 + i_2x_3 + i_1i_2x_4 = z_1 + i_2z_2$ and $\eta = y_1 + i_1y_2 + i_2y_3 + i_1i_2y_4 = w_1 + i_2w_2$ are two bicomplex numbers, then

$$\begin{aligned}
\overline{(\xi \times \eta)} &= \overline{\xi} \times \overline{\eta},\\
\widetilde{(\xi \times \eta)} &= \widetilde{\xi} \times \widetilde{\eta},\\
(\xi \times \eta)^\sharp &= \xi^\sharp \times \eta^\sharp.
\end{aligned}$$

**Definition 1.7.** Anjali et al., in [3], have defined two sets ${\mathbb{C}_2}^{m\times n}$ and ${\mathbb{C}_2}^n$ as follows:

$$\begin{aligned}
{\mathbb{C}_2}^{m\times n} &= \text{Set of all bicomplex matrices of order } m \times n.\\
&= \left\{A : A = [\xi_{ij}]_{m\times n}, \xi_{ij} \in \mathbb{C}_2, \forall i, j\right\}.\\
{\mathbb{C}_2}^n &= \{(\xi_1, \xi_2, \cdots, \xi_n) : \xi_1, \xi_2, \cdots, \xi_n \in \mathbb{C}_2\}.
\end{aligned}$$

Thus,

$$A = \begin{bmatrix} \xi_{11} & \xi_{12} & \cdots & \xi_{1n}\\ \xi_{21} & \xi_{22} & \cdots & \xi_{2n}\\ \vdots & \vdots & \ddots & \vdots\\ \xi_{m1} & \xi_{m2} & \cdots & \xi_{mn} \end{bmatrix}, \forall\ \xi_{ij} \in \mathbb{C}_2, i \in \{1, 2, \cdots, m\}, j \in \{1, 2, \cdots, n\}.$$

If $\xi_{ij} = {}^{ij}x_1 + i_1\,{}^{ij}x_2 + i_2\,{}^{ij}x_3 + i_1i_2\,{}^{ij}x_4 = {}^{ij}z_1 + i_2\,{}^{ij}z_2 = (\xi_{ij})^-\, e_1 + (\xi_{ij})^+\, e_2$, then due to the three forms of a bicomplex number, every bicomplex matrix $A = [\xi_{ij}]_{m\times n}$ can be uniquely expressed in the following three forms: the real form, the complex form, and the idempotent form, as follows:

$$\begin{aligned} A &= {}^{R}A_1 + i_1\,{}^{R}A_2 + i_2\,{}^{R}A_3 + i_1 i_2\,{}^{R}A_4 \rightarrow \text{ Real form of } A,\\ A &= {}^{C}A_1 + i_2\,{}^{C}A_2 \rightarrow \text{ Complex form of } A,\\ A &= A^{-}e_1 + A^{+}e_2 \rightarrow \text{ Idempotent form of } A. \end{aligned}$$

Where

$$\begin{gathered} {}^{R}A_1 = [{}^{ij}x_1]_{m\times n},\ {}^{R}A_2 = [{}^{ij}x_2]_{m\times n},\ {}^{R}A_3 = [{}^{ij}x_3]_{m\times n},\ {}^{R}A_4 = [{}^{ij}x_4]_{m\times n},\\ {}^{C}A_1 = {}^{R}A_1 + i_1\,{}^{R}A_2 = [{}^{ij}z_1]_{m\times n},\ {}^{C}A_2 = {}^{R}A_3 + i_1\,{}^{R}A_4 = [{}^{ij}z_2]_{m\times n},\\ A^{-} = [(\xi_{ij})^{-}]_{m\times n},\ A^{+} = [(\xi_{ij})^{+}]_{m\times n}. \end{gathered}$$

Let $\alpha \in \mathbb{C}_0$ or $\mathbb{C}_1$, and $P = (\xi_1, \xi_2, \cdots, \xi_n)$, $Q = (\eta_1, \eta_2, \cdots, \eta_n) \in {\mathbb{C}_2}^n$, and $A = [\xi_{ij}]_{m\times n}$, $B = [\eta_{ij}]_{m\times n} \in {\mathbb{C}_2}^{m\times n}$. Then the operations addition (+), scalar multiplication (·), and multiplication (×) on ${\mathbb{C}_2}^{m\times n}$ and ${\mathbb{C}_2}^n$ are defined as follows:

$$\begin{aligned} & A + B = [\xi_{ij} + \eta_{ij}]_{m\times n},\ \alpha\cdot A = [\alpha\cdot\xi_{ij}]_{m\times n},\\ & A\times B = \left({}^{C}A_1\,{}^{C}B_1 - {}^{C}A_2\,{}^{C}B_2\right) + i_2\left({}^{C}A_1\,{}^{C}B_2 + {}^{C}A_2\,{}^{C}B_1\right), \text{ provided that } m = n.\\ & P + Q = (\xi_1 + \eta_1, \xi_2 + \eta_2, \cdots, \xi_n + \eta_n), P\times Q = (\xi_1\times\eta_1, \xi_2\times\eta_2, \cdots, \xi_n\times\eta_n),\\ & \alpha P = (\alpha\xi_1, \alpha\xi_2, \cdots, \alpha\xi_n). \end{aligned}$$

**Theorem 1.8.** [3] *The structures* $({\mathbb{C}_2}^n, +, \cdot, \times, )$ *and* $\left({\mathbb{C}_2}^{n\times n}, +, \cdot, \times, \right)$ *form an algebra over* $\mathbb{C}_1$ *and* $\mathbb{C}_0$. *Moreover,* $\left({\mathbb{C}_2}^{m\times n}, +, \cdot\right)$ *forms a vector space over* $\mathbb{C}_1$ *and* $\mathbb{C}_0$.

## 2. Inner Products spaces on Bicomplex Numbers and Their Associated Spaces

This section expands the existing literature by introducing and analysing various inner product structures on $\mathbb{C}_2$ and its associated vector spaces. We develop the corresponding norm frameworks, explore their characteristics in comparison with the standard norm of a bicomplex number, and examine their roles in forming normed algebras and potential extensions toward bicomplex Hilbert spaces. Through this exploration, we aim to underscore the structural depth and analytical capabilities of bicomplex spaces, providing a foundation for continued research in this evolving field.

In the previous section, we discussed the three representations of a bicomplex number, namely the real representation, complex representation, and idempotent representation. Using the notions associated with these representations, we introduce three corresponding inner products on $\mathbb{C}_2$.

**Definition 2.1.** Let $\xi = x_1 + i_1x_2 + i_2x_3 + i_1i_2x_4 = z_1 + i_2z_2$ and $\eta = y_1 + i_1y_2 + i_2y_3 + i_1i_2y_4 = w_1 + i_2w_2$ be two arbitrary bicomplex numbers. Then, we define three maps $f_1, f_2$ and $f_3$ as follows:

$$f_1 : \mathbb{C}_2\times\mathbb{C}_2 \longrightarrow \mathbb{C}_1 \text{ such that } f_1(\xi,\eta) = \sum_{i=1}^{4} x_i y_i = x_1y_1 + x_2y_2 + x_3y_3 + x_4y_4, \tag{7}$$

$$f_2 : \mathbb{C}_2\times\mathbb{C}_2 \longrightarrow \mathbb{C}_1 \text{ such that } f_2(\xi,\eta) = z_1\overline{w_1} + z_2\overline{w_2}, \tag{8}$$

$$f_3 : \mathbb{C}_2\times\mathbb{C}_2 \longrightarrow \mathbb{C}_1 \text{ such that } f_3(\xi,\eta) = \xi^{-}\,(\overline{\eta^{-}}) + \xi^{+}\,(\overline{\eta^{+}}). \tag{9}$$

**Theorem 2.2.** *Let* $\xi = x_1 + i_1x_2 + i_2x_3 + i_1i_2x_4 = z_1 + i_2z_2$, $\eta = y_1 + i_1y_2 + i_2y_3 + i_1i_2y_4 = w_1 + i_2w_2$, *and* $\zeta = a_1 + i_1a_2 + i_2a_3 + i_1i_2a_4 = u_1 + i_2u_2$ *be the arbitrary bicomplex numbers, and* $\alpha = v_1 + i_1v_2 \in \mathbb{C}_1$. *Then, the map* $f_1$ *of the type* (7) *satisfies the following conditions:*

(a) $f_1(\xi,\xi) \geq 0$, *for all* $\xi \in \mathbb{C}_2$.
(b) $f_1(\xi,\xi) = 0$ *if and only if* $\xi = 0$.
(c) $\overline{f_1(\eta,\xi)} = f_1(\xi,\eta)$, *for all* $\xi,\eta \in \mathbb{C}_2$.
(d) $f_1(\xi+\eta,\zeta) = f_1(\xi,\zeta) + f_1(\eta,\zeta)$, *for all* $\xi,\eta,\zeta \in \mathbb{C}_2$.
(e) $f_1(\alpha\cdot\xi,\eta) = \alpha f_1(\xi,\eta)$, *for all* $\xi,\eta \in \mathbb{C}_2$ *and* $\alpha \in \mathbb{C}_1$.

*Proof.* For (a), we have

$$\begin{aligned} f_1(\xi,\xi) &= \sum_{i=1}^{4} x_i x_i \\ &= {x_1}^2 + {x_2}^2 + {x_3}^2 + {x_4}^2 \geq 0 \text{ as } x_1, x_2, x_3, x_4 \in \mathbb{C}_0. \end{aligned}$$

For (b), we put

$$\begin{aligned} & f_1(\xi,\xi) = 0 \\ \Longleftrightarrow\quad & \sum_{i=1}^{4} x_i x_i = 0 \\ \Longleftrightarrow\quad & {x_1}^2 + {x_2}^2 + {x_3}^2 + {x_4}^2 = 0 \\ \Longleftrightarrow\quad & x_1 = x_2 = x_3 = x_4 = 0, \text{ as } x_1, x_2, x_3, x_4 \in \mathbb{C}_0. \\ \Longleftrightarrow\quad & \xi = 0. \end{aligned}$$

For (c),

$$\begin{aligned} \overline{f_1(\eta,\xi)} &= \overline{\left\{ \sum_{i=1}^{4} y_i x_i \right\}} \\ &= \overline{y_1 x_1 + y_2 x_2 + y_3 x_3 + y_4 x_4} \\ &= y_1 x_1 + y_2 x_2 + y_3 x_3 + y_4 x_4, \text{ as } x_1, x_2, x_3, x_4, y_1, y_2, y_3, y_4 \in \mathbb{C}_0. \\ &= x_1 y_1 + x_2 y_2 + x_3 y_3 + x_4 y_4 \\ &= f_1(\xi,\eta). \end{aligned}$$

For (d),

$$\begin{aligned} f_1(\xi+\eta,\zeta) &= \sum_{i=1}^{4} (x_i + y_i) a_i \\ &= \sum_{i=1}^{4} x_i a_i + \sum_{i=1}^{4} y_i a_i \\ &= f_1(\xi,\zeta) + f_1(\eta,\zeta). \end{aligned}$$

For (e),

$$\begin{aligned} f_1(\alpha\cdot\xi,\eta) &= \sum_{i=1}^{4} (\alpha x_i) y_i \\ &= \alpha \sum_{i=1}^{4} x_i y_i \\ &= \alpha f_1(\xi,\eta). \end{aligned}$$

Hence, the theorem is proved. □

**Definition 2.3 (Real Inner Product on $\mathbb{C}_2$:).** Let $f_1$ be a map of the type (7). Then, by theorems 1.4 and 2.2, $f_1$ defines an inner product on $\mathbb{C}_2$. We refer to this inner product as the real inner product on $\mathbb{C}_2$ , and denote it by $\langle \bullet, \bullet \rangle_R$; that is,

$$f_1(\xi,\eta) = \langle \xi, \eta \rangle_R. \tag{10}$$

Consequently, $(\mathbb{C}_2, \langle \bullet, \bullet \rangle_R)$ is an inner product space. The norm on $\mathbb{C}_2$ induced by the real inner product is called the real norm and is denoted by $||\bullet||_R$; that is,

$$||\xi||_R = \sqrt{\langle \xi, \xi \rangle_R}. \tag{11}$$

**Remark 2.4.** Since the real norm of a bicomplex number $\xi$ is

$$\begin{aligned}||\xi||_R &= \sqrt{\langle \xi, \xi \rangle_R} \\ &= \left( \sum_{i=1}^{4} x_i^2 \right)^{\frac{1}{2}} \\ &= ||\xi||.\end{aligned}$$

Thus, the real norm coincides with the standard norm on $\mathbb{C}_2$.

**Theorem 2.5.** *Let $\xi = x_1 + i_1 x_2 + i_2 x_3 + i_1 i_2 x_4 = z_1 + i_2 z_2$, $\eta = y_1 + i_1 y_2 + i_2 y_3 + i_1 i_2 y_4 = w_1 + i_2 w_2$, and $\zeta = a_1 + i_1 a_2 + i_2 a_3 + i_1 i_2 a_4 = u_1 + i_2 u_2$ be the arbitrary bicomplex numbers, and $\alpha = v_1 + i_1 v_2 \in \mathbb{C}_1$. Then, the map $f_2$ of the type* (8) *satisfies the following conditions:*

*(a)* $f_2(\xi, \xi) \geq 0$, *for all* $\xi \in \mathbb{C}_2$.
*(b)* $f_2(\xi, \xi) = 0$ *if and only if* $\xi = 0$.
*(c)* $\overline{f_2(\eta, \xi)} = f_2(\xi, \eta)$, *for all* $\xi, \eta \in \mathbb{C}_2$.
*(d)* $f_2(\xi + \eta, \zeta) = f_2(\xi, \zeta) + f_2(\eta, \zeta)$, *for all* $\xi, \eta, \zeta \in \mathbb{C}_2$.
*(e)* $f_2(\alpha \cdot \xi, \eta) = \alpha f_2(\xi, \eta)$, *for all* $\xi, \eta \in \mathbb{C}_2$ *and* $\alpha \in \mathbb{C}_1$.

*Proof.* For (a), we have

$$\begin{aligned}f_2(\xi, \xi) &= z_1 \overline{z_1} + z_2 \overline{z_2} \\ &= |z_1|^2 + |z_2|^2 \geq 0.\end{aligned}$$

For (b), we put

$$\begin{aligned}& f_2(\xi, \xi) = 0 \\ \Longleftrightarrow\ & z_1 \overline{z_1} + z_2 \overline{z_2} = 0 \\ \Longleftrightarrow\ & |z_1|^2 + |z_2|^2 = 0 \\ \Longleftrightarrow\ & z_1 = z_2 = 0 \\ \Longleftrightarrow\ & \xi = 0.\end{aligned}$$

For (c),

$$\begin{aligned}\overline{f_2(\eta, \xi)} &= \overline{\{w_1 \overline{z_1} + w_2 \overline{z_2}\}} \\ &= z_1 \overline{w_1} + z_2 \overline{w_2} \\ &= f_2(\xi, \eta).\end{aligned}$$

For (d),

$$\begin{aligned}f_2(\xi + \eta, \zeta) &= (z_1 + w_1)\overline{u_1} + (z_2 + w_2)\overline{u_2} \\ &= z_1 \overline{u_1} + z_2 \overline{u_2} + w_1 \overline{u_1} + w_2 \overline{u_2} \\ &= f_2(\xi, \zeta) + f_2(\eta, \zeta).\end{aligned}$$

For (e),

$$\begin{aligned}f_2(\alpha \cdot \xi, \eta) &= (\alpha z_1)\overline{w_1} + (\alpha z_2)\overline{w_2} \\ &= \alpha(z_1 \overline{w_1} + z_2 \overline{w_2}) \\ &= \alpha f_2(\xi, \eta).\end{aligned}$$

Hence, the theorem is proved. □

**Definition 2.6 (Complex Inner Product on $\mathbb{C}_2$:).** Let $f_2$ be a map of the type (8). Then, by theorems 1.4 and 2.5, $f_2$ defines an inner product on $\mathbb{C}_2$. We refer to this inner product as the complex inner product on $\mathbb{C}_2$ , and denote it by $\langle \bullet, \bullet \rangle_C$; that is,

$$f_2(\xi, \eta) = \langle \xi, \eta \rangle_C. \tag{12}$$

Consequently, $(\mathbb{C}_2, \langle \bullet, \bullet \rangle_C)$ is an inner product space. The norm on $\mathbb{C}_2$ induced by the complex inner product is called the complex norm and is denoted by $||\bullet||_C$; that is,

$$||\xi||_C = \sqrt{\langle \xi, \xi \rangle_C}. \tag{13}$$

**Remark 2.7.** Since the complex norm of a bicomplex number $\xi$ is

$$\begin{aligned} ||\xi||_C &= \sqrt{\langle \xi, \xi \rangle_C} \\ &= \left(|z_1|^2 + |z_2|^2\right)^{\frac{1}{2}} \\ &= ||\xi||. \end{aligned}$$

Thus, the complex norm coincides with the standard norm on $\mathbb{C}_2$.

Theorem 2.8 follows directly from the result established in Theorems 2.4 and 2.7.

**Remark 2.8.** The real norm, the complex norm, and the standard norm on $\mathbb{C}_2$ coincide; that is,

$$||\xi|| = ||\xi||_R = ||\xi||_C.$$

**Theorem 2.9.** *Let $\xi = x_1 + i_1 x_2 + i_2 x_3 + i_1 i_2 x_4 = z_1 + i_2 z_2$, $\eta = y_1 + i_1 y_2 + i_2 y_3 + i_1 i_2 y_4 = w_1 + i_2 w_2$, and $\zeta = a_1 + i_1 a_2 + i_2 a_3 + i_1 i_2 a_4 = u_1 + i_2 u_2$ be the arbitrary bicomplex numbers, and $\alpha = v_1 + i_1 v_2 \in \mathbb{C}_1$. Then, the map $f_3$ of the type* (9) *satisfies the following conditions:*

*(a) $f_3(\xi, \xi) \geq 0$, for all $\xi \in \mathbb{C}_2$.*
*(b) $f_3(\xi, \xi) = 0$ if and only if $\xi = 0$.*
*(c) $\overline{f_3(\eta, \xi)} = f_3(\xi, \eta)$, for all $\xi, \eta \in \mathbb{C}_2$.*
*(d) $f_3(\xi + \eta, \zeta) = f_3(\xi, \zeta) + f_3(\eta, \zeta)$, for all $\xi, \eta, \zeta \in \mathbb{C}_2$.*
*(e) $f_3(\alpha \cdot \xi, \eta) = \alpha f_3(\xi, \eta)$, for all $\xi, \eta \in \mathbb{C}_2$ and $\alpha \in \mathbb{C}_1$.*

*Proof.* For (a), we have

$$\begin{aligned} f_3(\xi, \xi) &= \xi^- \left(\overline{\xi^-}\right) + \xi^+ \left(\overline{\xi^+}\right) \\ &= |\xi^-|^2 + |\xi^+|^2 \geq 0. \end{aligned}$$

For (b), we put

$$\begin{aligned} & f_3(\xi, \xi) = 0 \\ \iff & \xi^- \left(\overline{\xi^-}\right) + \xi^+ \left(\overline{\xi^+}\right) = 0 \\ \iff & |\xi^-|^2 + |\xi^+|^2 = 0 \\ \iff & \xi^- = \xi^+ = 0 \\ \iff & \xi = 0. \end{aligned}$$

For (c),

$$\begin{aligned} \overline{f_3(\eta, \xi)} &= \overline{\left\{\eta^- \left(\overline{\xi^-}\right) + \eta^+ \left(\overline{\xi^+}\right)\right\}} \\ &= \left(\overline{\eta^-}\right) \xi^- + \left(\overline{\eta^+}\right) \xi^+ \\ &= f_3(\xi, \eta). \end{aligned}$$

For (d),

$$\begin{aligned} f_3(\xi + \eta, \zeta) &= (\xi^- + \eta^-) \left(\overline{\zeta^-}\right) + (\xi^+ + \eta^+) \left(\overline{\zeta^+}\right) \\ &= \xi^- \left(\overline{\zeta^-}\right) + \xi^+ \left(\overline{\zeta^+}\right) + \eta^- \left(\overline{\zeta^-}\right) + \eta^+ \left(\overline{\zeta^+}\right) \\ &= f_3(\xi, \zeta) + f_3(\eta, \zeta). \end{aligned}$$

For (e),

$$\begin{aligned} f_3(\alpha \cdot \xi, \eta) &= (\alpha \xi^-) \left(\overline{\eta^-}\right) + (\alpha \xi^+) \left(\overline{\eta^+}\right) \\ &= \alpha \left\{\xi^- \left(\overline{\eta^-}\right) + \xi^+ \left(\overline{\eta^+}\right)\right\} \\ &= \alpha f_3(\xi, \eta). \end{aligned}$$

Hence, the theorem is proved. □

**Definition 2.10 (Idempotent Inner Product on $\mathbb{C}_2$:).** Let $f_3$ be a map of the type (9). Then, by theorems 1.4 and 2.9, $f_3$ defines an inner product on $\mathbb{C}_2$. We refer to this inner product as the idempotent inner product on $\mathbb{C}_2$ , and denote it by $\langle \bullet, \bullet \rangle_{ID}$; that is,

$$f_3(\xi, \eta) = \langle \xi, \eta \rangle_{ID}. \tag{14}$$

Consequently, $(\mathbb{C}_2, \langle \bullet, \bullet \rangle_{ID})$ is an inner product space. The norm on $\mathbb{C}_2$ induced by the idempotent inner product is called the idempotent norm and is denoted by $||\bullet||_{ID}$; that is,

$$||\xi||_{ID} = \sqrt{\langle \xi, \xi \rangle_{ID}}. \tag{15}$$

**Remark 2.11.** Since the idempotent norm of a bicomplex number $\xi$ is

$$\begin{aligned} ||\xi||_{ID} &= \sqrt{\langle \xi, \xi \rangle_{ID}} \\ &= \left(|\xi^-|^2 + |\xi^+|^2\right)^{\frac{1}{2}} \\ &= \sqrt{2}||\xi||. \end{aligned}$$

Thus, $||\xi||_{ID} = \sqrt{2}||\xi|| = \sqrt{2}||\xi||_R = \sqrt{2}||\xi||_C$, by theorem 2.8.

**Remark 2.12 (Compatibility of the Norms with Multiplication:).** Let $\xi$ and $\eta$ be two bicomplex numbers. Then, by theorem 1.3,

$$||\xi \times \eta|| \leq \sqrt{2}\ ||\xi||\ ||\eta|| \tag{16}$$

$$\Rightarrow ||\xi \times \eta||_R \leq \sqrt{2}\ ||\xi||_R\ ||\eta||_R, ||\xi \times \eta||_C \leq \sqrt{2}\ ||\xi||_C\ ||\eta||_C, \text{ by theorem 2.8.} \tag{17}$$

In theorem 1.4, we discussed the fact that the structure $(\mathbb{C}_2, +, \cdot, \times, ||\bullet||)$ forms a modified Banach algebra. Thus, by (17), the structures $(\mathbb{C}_2, +, \cdot, \times, ||\bullet||_R)$ and $(\mathbb{C}_2, +, \cdot, \times, ||\bullet||_C)$ form a modified Banach algebra, and the structures $(\mathbb{C}_2, +, \cdot, \langle \bullet, \bullet \rangle_R)$ and $(\mathbb{C}_2, +, \cdot, \langle \bullet, \bullet \rangle_C)$ are Hilbert spaces as well. Now

$$\begin{aligned} ||\xi \times \eta||_{ID} &= \sqrt{2}\ ||\xi \times \eta||, \text{ by theorem 2.11.} \\ \Rightarrow ||\xi \times \eta||_{ID} &\leq \sqrt{2}\ \sqrt{2}\ ||\xi||\ ||\eta||, \text{ by (16).} \\ \Rightarrow ||\xi \times \eta||_{ID} &\leq \sqrt{2}\ ||\xi||\ \sqrt{2}\ \ ||\eta|| \\ \Rightarrow ||\xi \times \eta||_{ID} &\leq ||\xi||_{ID}\ ||\eta||_{ID}, \text{ by theorem 2.11.} \end{aligned}$$

Let $\{\zeta_k\}_{k=1}^{\infty}$ be an arbitrary Cauchy sequence in $\mathbb{C}_2$ with respect to idempotent norm. Let $\epsilon > 0$, then there exists $n_0 \in \mathbb{N}$ such that

$$\begin{aligned} & \left|\left|\zeta_\alpha - \zeta_\beta\right|\right|_{ID} < \sqrt{2}\epsilon, \text{ for all } \alpha, \beta \geq n_0. \\ \Rightarrow\ & \sqrt{2}\ \left|\left|\zeta_\alpha - \zeta_\beta\right|\right| < \sqrt{2}\epsilon, \text{ for all } \alpha, \beta \geq n_0, \text{ by theorem 2.11.} \\ \Rightarrow\ & \left|\left|\zeta_\alpha - \zeta_\beta\right|\right| < \epsilon, \text{ for all } \alpha, \beta \geq n_0. \end{aligned}$$

Thus, $\{\zeta_k\}_{k=1}^{\infty}$ is the Cauchy sequence in $\mathbb{C}_2$ with respect to standard norm.

Now we consider an $\epsilon_1 > 0$. Then $\frac{\epsilon_1}{\sqrt{2}} > 0$ and by theorem 1.4, there exists $\zeta \in \mathbb{C}_2$ and $n_1 \in \mathbb{N}$ such that

$$\begin{aligned} & ||\zeta_\alpha - \zeta|| < \frac{\epsilon_1}{\sqrt{2}}, \text{ for all } \alpha \geq n_1. \\ \Rightarrow\ & \sqrt{2}\ ||\zeta_\alpha - \zeta|| < \epsilon_1, \text{ for all } \alpha \geq n_1. \\ \Rightarrow\ & ||\zeta_\alpha - \zeta||_{ID} < \epsilon_1, \text{ for all } \alpha \geq n_1, \text{ by theorem 2.11.} \\ \Rightarrow\ & ||\zeta_\alpha - \zeta||_{ID} \to 0, \text{ as } \alpha \to \infty. \\ \Rightarrow\ & \zeta_\alpha \to \zeta, \text{ as } \alpha \to \infty, \text{ with respect to idempotent norm.} \end{aligned}$$

This shows that every Cauchy sequence in $\mathbb{C}_2$ converges to an element of $\mathbb{C}_2$ under idempotent norm.

Hence, the structure $(\mathbb{C}_2, +, \cdot, \times, ||\bullet||_{ID})$ forms a Banach algebra and the structure $(\mathbb{C}_2, +, \cdot, \langle \bullet, \bullet \rangle_{ID})$ is a Hilbert space, whereas $(\mathbb{C}_2, +, \cdot, \times, ||\bullet||)$, $(\mathbb{C}_2, +, \cdot, \times, ||\bullet||_R)$, and $(\mathbb{C}_2, +, \cdot, \times, ||\bullet||_C)$ are modified Banach algebras. Thus, we have obtained a norm on $\mathbb{C}_2$ with respect to which $\mathbb{C}_2$ forms a Banach

algebra, thereby improving upon the previously known modified Banach algebra structure. This highlights a fundamental distinction between the idempotent norm and standard, real and complex norms.

**Remark 2.13.** Let $\xi = 1 + i_1 + 2i_2 + i_1 i_2$ and $\eta = 3 - 2i_1 + 3i_2 + 5i_1 i_2$. Then

$$\begin{aligned}
\langle \xi, \eta \rangle_R &= 12, \\
\langle \xi, \eta \rangle_C &= 12 - 2i_1, \\
\langle \xi, \eta \rangle_{ID} &= 24 - 4i_1.
\end{aligned}$$

This shows that the real, complex, and idempotent inner products on $\mathbb{C}_2$ are pairwise distinct.

Using the notions associated with the representations of a bicomplex number, we introduce three corresponding inner products on ${\mathbb{C}_2}^n$.

**Definition 2.14.** Let $P = (\xi_1, \xi_2, \cdots, \xi_n)$ and $Q = (\eta_1, \eta_2, \cdots, \eta_n)$ be two arbitrary elements of ${\mathbb{C}_2}^n$ such that

$$\begin{aligned}
\xi_k &= {}^k x_1 + i_1\, {}^k x_2 + i_2\, {}^k x_3 + i_1 i_2\, {}^k x_4 = {}^k z_1 + i_2\, {}^k z_2 = (\xi_k)^-\, e_1 + (\xi_k)^+\, e_2, \\
\eta_k &= {}^k y_1 + i_1\, {}^k y_2 + i_2\, {}^k y_3 + i_1 i_2\, {}^k y_4 = {}^k w_1 + i_2\, {}^k w_2 = (\eta_k)^-\, e_1 + (\eta_k)^+\, e_2, k = 1, 2, \cdots, n.
\end{aligned}$$

Then, we define three maps $f_4, f_5$ and $f_6$ as follows:

$$f_4, f_5, f_6 : {\mathbb{C}_2}^n \times {\mathbb{C}_2}^n \longrightarrow \mathbb{C}_1 \text{ such that}$$

$$f_4(P, Q) = \sum_{k=1}^{n} ({}^k x_1)({}^k y_1) + ({}^k x_2)({}^k y_2) + ({}^k x_3)({}^k y_3) + ({}^k x_4)({}^k y_4), \tag{18}$$

$$f_5(P, Q) = \sum_{k=1}^{n} ({}^k z_1)\overline{({}^k w_1)} + ({}^k z_2)\overline{({}^k w_2)}, \tag{19}$$

$$f_6(P, Q) = \sum_{k=1}^{n} (\xi_k)^-\, \overline{(\eta_k)^-} + (\xi_k)^+\, \overline{(\eta_k)^+}. \tag{20}$$

**Theorem 2.15.** *Let $P = (\xi_1, \xi_2, \cdots, \xi_n)$, $Q = (\eta_1, \eta_2, \cdots, \eta_n)$, and $R = (\zeta_1, \zeta_2, \cdots, \zeta_n)$ be the arbitrary elements of ${\mathbb{C}_2}^n$ such that*

$$\begin{aligned}
\xi_k &= {}^k x_1 + i_1\, {}^k x_2 + i_2\, {}^k x_3 + i_1 i_2\, {}^k x_4 = {}^k z_1 + i_2\, {}^k z_2 = (\xi_k)^-\, e_1 + (\xi_k)^+\, e_2, \\
\eta_k &= {}^k y_1 + i_1\, {}^k y_2 + i_2\, {}^k y_3 + i_1 i_2\, {}^k y_4 = {}^k w_1 + i_2\, {}^k w_2 = (\eta_k)^-\, e_1 + (\eta_k)^+\, e_2, \\
\zeta_k &= {}^k a_1 + i_1\, {}^k a_2 + i_2\, {}^k a_3 + i_1 i_2\, {}^k a_4 = {}^k u_1 + i_2\, {}^k u_2 = (\zeta_k)^-\, e_1 + (\zeta_k)^+\, e_2, k = 1, 2, \cdots, n.
\end{aligned}$$

*and $\alpha = v_1 + i_1 v_2 \in \mathbb{C}_1$. Then, the map $f_4$ of the type* (18) *satisfies the following conditions:*

*(a)* $f_4(P, P) \geq 0$, *for all* $P \in {\mathbb{C}_2}^n$.
*(b)* $f_4(P, P) = 0$ *if and only if* $P = 0$.
*(c)* $\overline{f_4(Q, P)} = f_4(P, Q)$, *for all* $P, Q \in {\mathbb{C}_2}^n$.
*(d)* $f_4(P + Q, R) = f_4(P, R) + f_4(Q, R)$, *for all* $P, Q, R \in {\mathbb{C}_2}^n$.
*(e)* $f_4(\alpha P, Q) = \alpha f_4(P, Q)$, *for all* $P, Q \in {\mathbb{C}_2}^n$ *and* $\alpha \in \mathbb{C}_1$.

*Proof.* Since

$$\begin{aligned}
f_4(P, Q) &= \sum_{k=1}^{n} ({}^k x_1)({}^k y_1) + ({}^k x_2)({}^k y_2) + ({}^k x_3)({}^k y_3) + ({}^k x_4)({}^k y_4) \\
\Rightarrow f_4(P, Q) &= \sum_{k=1}^{n} \langle \xi_k, \eta_k \rangle_R, \text{ by (7) and (10).}
\end{aligned}$$

For (a), we have

$$f_4(P, P) = \sum_{k=1}^{n} \langle \xi_k, \xi_k \rangle_R.$$

Since $\langle \xi_k, \xi_k \rangle_R \geq 0$, for all $k \in \{1, 2, \cdots, n\}$. Therefore

$$f_4(P, P) \;=\; \sum_{k=1}^{n} \langle \xi_k, \xi_k \rangle_R \geq 0.$$

For (b), we put

$$\begin{aligned}
& f_4(P, P) = 0 \\
\iff & \sum_{k=1}^{n} \langle \xi_k, \xi_k \rangle_R = 0 \\
\iff & \langle \xi_k, \xi_k \rangle_R = 0, \text{ for all } k \in \{1, 2, \cdots, n\}. \\
\iff & \xi_k = 0, \text{ for all } k \in \{1, 2, \cdots, n\}. \\
\iff & P = (\xi_1, \xi_2, \cdots, \xi_n) = \{0, 0, \cdots, 0\}.
\end{aligned}$$

For (c),

$$\begin{aligned}
\overline{f_4(Q, P)} &= \overline{\sum_{k=1}^{n} \langle \eta_k, \xi_k \rangle_R} \\
&= \sum_{k=1}^{n} \overline{\langle \eta_k, \xi_k \rangle_R} \\
&= \sum_{k=1}^{n} \langle \xi_k, \eta_k \rangle_R \\
&= f_4(P, Q).
\end{aligned}$$

For (d),

$$\begin{aligned}
f_4(P + Q, R) &= \sum_{k=1}^{n} \langle \xi_k + \eta_k, \zeta_k \rangle_R \\
&= \sum_{k=1}^{n} \left\{ \langle \xi_k, \zeta_k \rangle_R + \langle \eta_k, \zeta_k \rangle_R \right\} \\
&= \sum_{k=1}^{n} \langle \xi_k, \zeta_k \rangle_R + \sum_{k=1}^{n} \langle \eta_k, \zeta_k \rangle_R \\
&= f_4(P, R) + f_4(Q, R).
\end{aligned}$$

For (e),

$$\begin{aligned}
f_4(\alpha P, Q) &= \sum_{k=1}^{n} \langle \alpha \xi_k, \eta_k \rangle_R \\
&= \sum_{k=1}^{n} \alpha \langle \xi_k, \eta_k \rangle_R \\
&= \alpha \sum_{k=1}^{n} \langle \xi_k, \eta_k \rangle_R \\
&= \alpha f_4(P, Q).
\end{aligned}$$

Hence, the theorem is proved. □

**Definition 2.16** (**Real Inner Product on** $\mathbb{C}_2{}^n$**:**)**.** Let $f_4$ be a map of the type (18). Then, by theorems 1.8 and 2.15, $f_4$ defines an inner product on $\mathbb{C}_2{}^n$. We refer to this inner product as the real inner product on $\mathbb{C}_2{}^n$ , and denote it by $\langle\langle \bullet, \bullet \rangle\rangle_R$; that is,

$$f_4(P, Q) = \langle\langle P, Q \rangle\rangle_R.$$

Consequently, $(\mathbb{C}_2{}^n, \langle\langle \bullet, \bullet \rangle\rangle_R)$ is an inner product space. The norm on $\mathbb{C}_2{}^n$ induced by the real inner product is called the real norm and is denoted $|||\bullet|||_R$; that is,

$$|||P|||_R = \sqrt{\langle\langle P, P \rangle\rangle_R}.$$

**Remark 2.17.** Since the real norm of $P$ is

$$\begin{aligned} |||P|||_R &= \sqrt{\langle\langle P, P \rangle\rangle_R} \\ &= \left( \sum_{k=1}^{n} \langle \xi_k, \xi_k \rangle_R \right)^{\frac{1}{2}} \\ &= \left( \sum_{k=1}^{n} (||\xi_k||_R)^2 \right)^{\frac{1}{2}}, \text{ by (11).} \\ &= \left( \sum_{k=1}^{n} ||\xi_k||^2 \right)^{\frac{1}{2}}, \text{ by theorem 2.4.} \end{aligned}$$

Thus, $|||P|||_R = \left(||\xi_1||^2 + ||\xi_2||^2 + \cdots + ||\xi_n||^2\right)^{\frac{1}{2}}$.

**Theorem 2.18.** *Let $P = (\xi_1, \xi_2, \cdots, \xi_n)$, $Q = (\eta_1, \eta_2, \cdots, \eta_n)$, and $R = (\zeta_1, \zeta_2, \cdots, \zeta_n)$ be the arbitrary elements of $\mathbb{C}_2{}^n$ such that*

$$\begin{aligned} \xi_k &= {}^k x_1 + i_1\, {}^k x_2 + i_2\, {}^k x_3 + i_1 i_2\, {}^k x_4 = {}^k z_1 + i_2\, {}^k z_2 = (\xi_k)^-\, e_1 + (\xi_k)^+\, e_2, \\ \eta_k &= {}^k y_1 + i_1\, {}^k y_2 + i_2\, {}^k y_3 + i_1 i_2\, {}^k y_4 = {}^k w_1 + i_2\, {}^k w_2 = (\eta_k)^-\, e_1 + (\eta_k)^+\, e_2, \\ \zeta_k &= {}^k a_1 + i_1\, {}^k a_2 + i_2\, {}^k a_3 + i_1 i_2\, {}^k a_4 = {}^k u_1 + i_2\, {}^k u_2 = (\zeta_k)^-\, e_1 + (\zeta_k)^+\, e_2, k = 1, 2, \cdots, n. \end{aligned}$$

*and $\alpha = v_1 + i_1 v_2 \in \mathbb{C}_1$. Then, the map $f_5$ of the type* (19) *satisfies the following conditions:*

*(a) $f_5(P, P) \geq 0$, for all $P \in \mathbb{C}_2{}^n$.*
*(b) $f_5(P, P) = 0$ if and only if $P = 0$.*
*(c) $\overline{f_5(Q, P)} = f_5(P, Q)$, for all $P, Q \in \mathbb{C}_2{}^n$.*
*(d) $f_5(P + Q, R) = f_5(P, R) + f_5(Q, R)$, for all $P, Q, R \in \mathbb{C}_2{}^n$.*
*(e) $f_5(\alpha P, Q) = \alpha f_5(P, Q)$, for all $P, Q \in \mathbb{C}_2{}^n$ and $\alpha \in \mathbb{C}_1$.*

*Proof.* Since

$$\begin{aligned} f_5(P, Q) &= \sum_{k=1}^{n} ({}^k z_1)\overline{({}^k w_1)} + ({}^k z_2)\overline{({}^k w_2)} \\ \Rightarrow f_5(P, Q) &= \sum_{k=1}^{n} \langle \xi_k, \eta_k \rangle_C, \text{ by (8) and (12).} \end{aligned}$$

For (a), we have

$$f_5(P, P) = \sum_{k=1}^{n} \langle \xi_k, \xi_k \rangle_C.$$

Since $\langle \xi_k, \xi_k \rangle_C \geq 0$, for all $k \in \{1, 2, \cdots, n\}$. Therefore

$$f_5(P, P) = \sum_{k=1}^{n} \langle \xi_k, \xi_k \rangle_C \geq 0.$$

For (b), we put

$$\begin{aligned} & f_5(P,P)=0 \\ \Longleftrightarrow\quad & \sum_{k=1}^{n}\langle\xi_k,\xi_k\rangle_C=0 \\ \Longleftrightarrow\quad & \langle\xi_k,\xi_k\rangle_C=0,\ \text{for all } k\in\{1,2,\cdots,n\}. \\ \Longleftrightarrow\quad & \xi_k=0,\ \text{for all } k\in\{1,2,\cdots,n\}. \\ \Longleftrightarrow\quad & P=(\xi_1,\xi_2,\cdots,\xi_n)=\{0,0,\cdots,0\}. \end{aligned}$$

For (c),

$$\begin{aligned} \overline{f_5(Q,P)} &= \overline{\sum_{k=1}^{n}\langle\eta_k,\xi_k\rangle_C} \\ &= \sum_{k=1}^{n}\overline{\langle\eta_k,\xi_k\rangle_C} \\ &= \sum_{k=1}^{n}\langle\xi_k,\eta_k\rangle_C \\ &= f_5(P,Q). \end{aligned}$$

For (d),

$$\begin{aligned} f_5(P+Q,R) &= \sum_{k=1}^{n}\langle\xi_k+\eta_k,\zeta_k\rangle_C \\ &= \sum_{k=1}^{n}\left\{\langle\xi_k,\zeta_k\rangle_C+\langle\eta_k,\zeta_k\rangle_C\right\} \\ &= \sum_{k=1}^{n}\langle\xi_k,\zeta_k\rangle_C+\sum_{k=1}^{n}\langle\eta_k,\zeta_k\rangle_C \\ &= f_5(P,R)+f_5(Q,R). \end{aligned}$$

For (e),

$$\begin{aligned} f_5(\alpha P,Q) &= \sum_{k=1}^{n}\langle\alpha\xi_k,\eta_k\rangle_C \\ &= \sum_{k=1}^{n}\alpha\langle\xi_k,\eta_k\rangle_C \\ &= \alpha\sum_{k=1}^{n}\langle\xi_k,\eta_k\rangle_C \\ &= \alpha f_5(P,Q). \end{aligned}$$

Hence, the theorem is proved. □

**Definition 2.19 (Complex Inner Product on $\mathbb{C}_2{}^n$:).** Let $f_5$ be a map of the type (19). Then, by theorems 1.8 and 2.18, $f_5$ defines an inner product on $\mathbb{C}_2{}^n$. We refer to this inner product as the complex inner product on $\mathbb{C}_2{}^n$ , and denote it by $\langle\langle\bullet,\bullet\rangle\rangle_C$; that is,

$$f_5(P,Q)=\langle\langle P,Q\rangle\rangle_C.$$

Consequently, $(\mathbb{C}_2{}^n,\langle\langle\bullet,\bullet\rangle\rangle_C)$ is an inner product space. The norm on $\mathbb{C}_2{}^n$ induced by the complex inner product is called the complex norm and is denoted $|||\bullet|||_C$; that is,

$$|||P|||_C=\sqrt{\langle\langle P,P\rangle\rangle_C}.$$

**Remark 2.20.** Since the complex norm of $P$ is

$$\begin{aligned} |||P|||_C &= \sqrt{\langle\langle P,P\rangle\rangle_C} \\ &= \left(\sum_{k=1}^{n}\langle\xi_k,\xi_k\rangle_C\right)^{\frac{1}{2}} \\ &= \left(\sum_{k=1}^{n}(||\xi_k||_C)^2\right)^{\frac{1}{2}}, \text{ by (13).} \\ &= \left(\sum_{k=1}^{n}||\xi_k||^2\right)^{\frac{1}{2}}, \text{ by theorem 2.7.} \end{aligned}$$

Thus, $|||P|||_C = \left(||\xi_1||^2 + ||\xi_2||^2 + \cdots + ||\xi_n||^2\right)^{\frac{1}{2}} = |||P|||_R$.

**Theorem 2.21.** *Let $P = (\xi_1, \xi_2, \cdots, \xi_n)$, $Q = (\eta_1, \eta_2, \cdots, \eta_n)$, and $R = (\zeta_1, \zeta_2, \cdots, \zeta_n)$ be the arbitrary elements of ${\mathbb{C}_2}^n$ such that*

$$\begin{aligned} \xi_k &= {}^kx_1 + i_1\,{}^kx_2 + i_2\,{}^kx_3 + i_1i_2\,{}^kx_4 = {}^kz_1 + i_2\,{}^kz_2 = (\xi_k)^-\,e_1 + (\xi_k)^+\,e_2, \\ \eta_k &= {}^ky_1 + i_1\,{}^ky_2 + i_2\,{}^ky_3 + i_1i_2\,{}^ky_4 = {}^kw_1 + i_2\,{}^kw_2 = (\eta_k)^-\,e_1 + (\eta_k)^+\,e_2, \\ \zeta_k &= {}^ka_1 + i_1\,{}^ka_2 + i_2\,{}^ka_3 + i_1i_2\,{}^ka_4 = {}^ku_1 + i_2\,{}^ku_2 = (\zeta_k)^-\,e_1 + (\zeta_k)^+\,e_2, k = 1, 2, \cdots, n. \end{aligned}$$

*and $\alpha = v_1 + i_1v_2 \in \mathbb{C}_1$. Then, the map $f_6$ of the type* (20) *satisfies the following conditions:*

*(a)* $f_6(P,P) \geq 0$, *for all* $P \in {\mathbb{C}_2}^n$.
*(b)* $f_6(P,P) = 0$ *if and only if* $P = 0$.
*(c)* $\overline{f_6(Q,P)} = f_6(P,Q)$, *for all* $P,Q \in {\mathbb{C}_2}^n$.
*(d)* $f_6(P+Q,R) = f_6(P,R) + f_6(Q,R)$, *for all* $P,Q,R \in {\mathbb{C}_2}^n$.
*(e)* $f_6(\alpha P,Q) = \alpha f_6(P,Q)$, *for all* $P,Q \in {\mathbb{C}_2}^n$ *and* $\alpha \in \mathbb{C}_1$.

*Proof.* Since

$$\begin{aligned} f_6(P,Q) &= \sum_{k=1}^{n}(\xi_k)^-\,\overline{(\eta_k)^-} + (\xi_k)^+\,\overline{(\eta_k)^+} \\ \Rightarrow f_6(P,Q) &= \sum_{k=1}^{n}\langle\xi_k,\eta_k\rangle_{ID}, \text{ by (9) and (14).} \end{aligned}$$

For (a), we have

$$f_6(P,P) = \sum_{k=1}^{n}\langle\xi_k,\xi_k\rangle_{ID}.$$

Since $\langle\xi_k,\xi_k\rangle_{ID} \geq 0$, for all $k \in \{1, 2, \cdots, n\}$. Therefore

$$f_6(P,P) = \sum_{k=1}^{n}\langle\xi_k,\xi_k\rangle_{ID} \geq 0.$$

For (b), we put

$$\begin{aligned} & f_6(P,P) = 0 \\ \Longleftrightarrow\ & \sum_{k=1}^{n}\langle\xi_k,\xi_k\rangle_{ID} = 0 \\ \Longleftrightarrow\ & \langle\xi_k,\xi_k\rangle_{ID} = 0, \text{ for all } k \in \{1, 2, \cdots, n\}. \\ \Longleftrightarrow\ & \xi_k = 0, \text{ for all } k \in \{1, 2, \cdots, n\}. \\ \Longleftrightarrow\ & P = (\xi_1, \xi_2, \cdots, \xi_n) = \{0, 0, \cdots, 0\}. \end{aligned}$$

For (c),

$$\begin{aligned}
\overline{f_6(Q,P)} &= \overline{\sum_{k=1}^{n}\langle \eta_k, \xi_k\rangle_{ID}} \\
&= \sum_{k=1}^{n}\overline{\langle \eta_k, \xi_k\rangle_{ID}} \\
&= \sum_{k=1}^{n}\langle \xi_k, \eta_k\rangle_{ID} \\
&= f_6(P,Q).
\end{aligned}$$

For (d),

$$\begin{aligned}
f_6(P+Q,R) &= \sum_{k=1}^{n}\langle \xi_k+\eta_k, \zeta_k\rangle_{ID} \\
&= \sum_{k=1}^{n}\left\{\langle \xi_k, \zeta_k\rangle_{ID}+\langle \eta_k, \zeta_k\rangle_{ID}\right\} \\
&= \sum_{k=1}^{n}\langle \xi_k, \zeta_k\rangle_{ID}+\sum_{k=1}^{n}\langle \eta_k, \zeta_k\rangle_{ID} \\
&= f_6(P,R)+f_6(Q,R).
\end{aligned}$$

For (e),

$$\begin{aligned}
f_6(\alpha P,Q) &= \sum_{k=1}^{n}\langle \alpha\xi_k, \eta_k\rangle_{ID} \\
&= \sum_{k=1}^{n}\alpha\langle \xi_k, \eta_k\rangle_{ID} \\
&= \alpha\sum_{k=1}^{n}\langle \xi_k, \eta_k\rangle_{ID} \\
&= \alpha f_6(P,Q).
\end{aligned}$$

Hence, the theorem is proved. □

**Definition 2.22 (Idempotent Inner Product on $\mathbb{C}_2{}^n$:).** Let $f_6$ be a map of the type (20). Then, by theorems 1.8 and 2.21, $f_6$ defines an inner product on $\mathbb{C}_2{}^n$. We refer to this inner product as the idempotent inner product on $\mathbb{C}_2{}^n$ , and denote it by $\langle\langle \bullet, \bullet\rangle\rangle_{ID}$; that is,

$$f_6(P,Q) = \langle\langle P,Q\rangle\rangle_{ID}.$$

Consequently, $(\mathbb{C}_2{}^n, \langle\langle \bullet, \bullet\rangle\rangle_{ID})$ is an inner product space. The norm on $\mathbb{C}_2{}^n$ induced by the idempotent inner product is called the idempotent norm and is denoted $|||\bullet|||_{ID}$; that is,

$$|||P|||_{ID} = \sqrt{\langle\langle P,P\rangle\rangle_{ID}}.$$

**Remark 2.23.** Since the idempotent norm of $P$ is

$$\begin{aligned} |||P|||_{ID} &= \sqrt{\langle\langle P,P\rangle\rangle_{ID}} \\ &= \left(\sum_{k=1}^{n}\langle\xi_k,\xi_k\rangle_{ID}\right)^{\frac{1}{2}} \\ &= \left(\sum_{k=1}^{n}(||\xi_k||_{ID})^2\right)^{\frac{1}{2}}, \text{ by (15).} \\ &= \sqrt{2}\left(\sum_{k=1}^{n}||\xi_k||^2\right)^{\frac{1}{2}}, \text{ by theorem 2.11.} \end{aligned}$$

Thus, $|||P|||_{ID} = \sqrt{2}\big(||\xi_1||^2 + ||\xi_2||^2 + \cdots + ||\xi_n||^2\big)^{\frac{1}{2}} = \sqrt{2}|||P|||_R = \sqrt{2}|||P|||_C$, by theorem 2.20.

**Remark 2.24.** Let $P = (\xi_1, \xi_2)$ and $Q = (\eta_1, \eta_2)$ such that

$$\begin{aligned} \xi_1 &= 1 + i_1 + 2i_2 + i_1i_2, \\ \xi_2 &= 2 + i_1 + i_2 + 3i_1i_2, \\ \eta_1 &= 3 - 2i_1 + 3i_2 + 5i_1i_2, \\ \eta_2 &= 0 + 3i_1 + i_2 - 3i_1i_2. \end{aligned}$$

Then

$$\begin{aligned} \langle\langle P,Q\rangle\rangle_R &= \langle\xi_1,\eta_1\rangle_R + \langle\xi_2,\eta_2\rangle_R = 12 - 5 = 7, \\ \langle\langle P,Q\rangle\rangle_C &= \langle\xi_1,\eta_1\rangle_C + \langle\xi_2,\eta_2\rangle_C = 12 - 2i_1 - 5 = 7 - 2i_1, \\ \langle\langle P,Q\rangle\rangle_{ID} &= \langle\xi_1,\eta_1\rangle_{ID} + \langle\xi_2,\eta_2\rangle_{ID} = 24 - 4i_1 - 10 = 14 - 4i_1. \end{aligned}$$

This shows that the real, complex, and idempotent inner products on ${\mathbb{C}_2}^n$ are pairwise distinct.

**Theorem 2.25.** *Let $P = (\xi_1, \xi_2, \cdots, \xi_n)$ and $Q = (\eta_1, \eta_2, \cdots, \eta_n)$ be two arbitrary elements of ${\mathbb{C}_2}^n$. Then*

*(a)* $|||PQ|||_R \leq \sqrt{2}\ |||P|||_R\ |||Q|||_R$.
*(b)* $|||PQ|||_C \leq \sqrt{2}\ |||P|||_C\ |||Q|||_C$.
*(c)* $|||PQ|||_{ID} \leq |||P|||_{ID}\ |||Q|||_{ID}$.

*Moreover, the constant $\sqrt{2}$ in the inequalities (a) and (b) is best possible.*

*Proof.* Since $|||P|||_R = \big(||\xi_1||^2 + ||\xi_2||^2 + \cdots + ||\xi_n||^2\big)^{\frac{1}{2}} = |||P|||_C$, by theorem 2.20. Therefore

$$\begin{aligned} (|||PQ|||_R)^2 &= \big(||\xi_1\eta_1||^2 + ||\xi_2\eta_2||^2 + \cdots + ||\xi_n\eta_n||^2\big) \\ \Rightarrow (|||PQ|||_R)^2 &\leq \big(2||\xi_1||^2||\eta_1||^2 + 2||\xi_2||^2||\eta_2||^2 + \cdots + 2||\xi_n||^2||\eta_n||^2\big), \text{ by theorem 1.3.} \\ \Rightarrow (|||PQ|||_R)^2 &\leq 2\big(||\xi_1||^2||\eta_1||^2 + ||\xi_2||^2||\eta_2||^2 + \cdots + ||\xi_n||^2||\eta_n||^2\big) \\ \Rightarrow (|||PQ|||_R)^2 &\leq 2\big(||\xi_1||^2 + ||\xi_2||^2 + \cdots + ||\xi_n||^2\big)\big(||\eta_1||^2 + ||\eta_2||^2 + \cdots + ||\eta_n||^2\big) \\ \Rightarrow (|||PQ|||_R)^2 &\leq 2(|||P|||_R)^2\ (|||Q|||_R)^2. \end{aligned}$$

Therefore

$$|||PQ|||_R \leq \sqrt{2}\ |||P|||_R\ |||Q|||_R. \tag{21}$$

In the same manner, since $|||P|||_{ID} = \sqrt{2}|||P|||_R = \sqrt{2}|||P|||_C$, by theorem 2.23. Therefore

$$\begin{aligned} |||PQ|||_{ID} &= \sqrt{2}|||PQ|||_R \\ \Rightarrow |||PQ|||_{ID} &\leq \sqrt{2}\sqrt{2}\ |||P|||_R\ |||Q|||_R, \text{ by (21).} \\ \Rightarrow |||PQ|||_{ID} &\leq \sqrt{2}\ |||P|||_R\ \sqrt{2}\ |||Q|||_R \\ \Rightarrow |||PQ|||_{ID} &\leq |||P|||_{ID}\ |||Q|||_{ID} \end{aligned}$$

Moreover, if $P = (e_1, 0, \cdots, 0)$ and $Q = (e_1, 0, \cdots, 0)$ are the elements of ${\mathbb{C}_2}^n$, then

$$\begin{aligned} |||PQ|||_R &= \left(||e_1||^2 + ||0||^2 + \cdots + ||0||^2\right)^{\frac{1}{2}} \\ &= \frac{1}{\sqrt{2}} \\ &= \sqrt{2}\,|||P|||_R\,|||Q|||_R. \end{aligned}$$

Hence, the theorem is proved. □

**Theorem 2.26.** *Every Cauchy sequence in ${\mathbb{C}_2}^n$ converges to an element of ${\mathbb{C}_2}^n$ under real norm.*

*Proof.* Let $\{P_k\}_{k=1}^{\infty}$ be an arbitrary Cauchy sequence in ${\mathbb{C}_2}^n$ with respect to real norm, where $P_k = \left(\xi_{1(k)}, \xi_{2(k)}, \cdots, \xi_{n(k)}\right)$. Let $\epsilon > 0$, then there exists $n_0 \in \mathbb{N}$ such that

$$\left|\left|\left|P_\alpha - P_\beta\right|\right|\right|_R < \epsilon, \text{ for all natural numbers } \alpha \text{ and } \beta \text{ such that } \alpha, \beta \geq n_0.$$

Therefore, by theorem 2.20

$$\begin{aligned} &\left(\left|\left|\xi_{1(\alpha)} - \xi_{1(\beta)}\right|\right|^2 + \left|\left|\xi_{2(\alpha)} - \xi_{2(\beta)}\right|\right|^2 + \cdots + \left|\left|\xi_{n(\alpha)} - \xi_{n(\beta)}\right|\right|^2\right)^{\frac{1}{2}} < \epsilon, \text{ for all } \alpha, \beta \geq n_0. \\ \Rightarrow\ &\left(\left|\left|\xi_{1(\alpha)} - \xi_{1(\beta)}\right|\right|^2 + \left|\left|\xi_{2(\alpha)} - \xi_{2(\beta)}\right|\right|^2 + \cdots + \left|\left|\xi_{n(\alpha)} - \xi_{n(\beta)}\right|\right|^2\right) < \epsilon^2, \text{ for all } \alpha, \beta \geq n_0. \\ \Rightarrow\ &\left|\left|\xi_{1(\alpha)} - \xi_{1(\beta)}\right|\right|^2 < \epsilon^2,\ \left|\left|\xi_{2(\alpha)} - \xi_{2(\beta)}\right|\right|^2 < \epsilon^2, \cdots, \left|\left|\xi_{n(\alpha)} - \xi_{n(\beta)}\right|\right|^2 < \epsilon^2, \text{ for all } \alpha, \beta \geq n_0. \\ \Rightarrow\ &\left|\left|\xi_{1(\alpha)} - \xi_{1(\beta)}\right|\right| < \epsilon,\ \left|\left|\xi_{2(\alpha)} - \xi_{2(\beta)}\right|\right| < \epsilon, \cdots, \left|\left|\xi_{n(\alpha)} - \xi_{n(\beta)}\right|\right| < \epsilon, \text{ for all } \alpha, \beta \geq n_0. \end{aligned}$$

Thus, $\left\{\xi_{1(k)}\right\}_{k=1}^{\infty}, \left\{\xi_{2(k)}\right\}_{k=1}^{\infty}, \cdots, \left\{\xi_{n(k)}\right\}_{k=1}^{\infty}$ are the Cauchy sequences in $\mathbb{C}_2$ with respect to standard norm.

Now we consider an $\epsilon_1 > 0$. Then $\frac{\epsilon_1}{\sqrt{n}} > 0$ and by theorem 2.12, there exists $\xi_1, \xi_2, \cdots, \xi_n \in \mathbb{C}_2$ and $\gamma_1, \gamma_2, \cdots, \gamma_n \in \mathbb{N}$ such that

$$\begin{aligned} &\left|\left|\xi_{1(\alpha_1)} - \xi_1\right|\right| < \frac{\epsilon_1}{\sqrt{n}},\ \left|\left|\xi_{2(\alpha_2)} - \xi_2\right|\right| < \frac{\epsilon_1}{\sqrt{n}}, \cdots, \left|\left|\xi_{n(\alpha_n)} - \xi_n\right|\right| < \frac{\epsilon_1}{\sqrt{n}}, \text{ for all } \alpha_1 \geq \gamma_1, \\ &\alpha_2 \geq \gamma_2, \cdots, \alpha_n \geq \gamma_n. \\ \Rightarrow\ &\left|\left|\xi_{1(\alpha)} - \xi_1\right|\right| < \frac{\epsilon_1}{\sqrt{n}},\ \left|\left|\xi_{2(\alpha)} - \xi_2\right|\right| < \frac{\epsilon_1}{\sqrt{n}}, \cdots, \left|\left|\xi_{n(\alpha)} - \xi_n\right|\right| < \frac{\epsilon_1}{\sqrt{n}}, \text{ for all } \alpha \geq \gamma, \\ &\text{where } \gamma = \max\{\gamma_1, \gamma_2, \cdots, \gamma_n\}. \\ \Rightarrow\ &\left|\left|\xi_{1(\alpha)} - \xi_1\right|\right|^2 < \frac{{\epsilon_1}^2}{n}, \left|\left|\xi_{2(\alpha)} - \xi_2\right|\right|^2 < \frac{{\epsilon_1}^2}{n}, \cdots, \left|\left|\xi_{n(\alpha)} - \xi_n\right|\right|^2 < \frac{{\epsilon_1}^2}{n}, \text{ for all } \alpha \geq \gamma. \\ \Rightarrow\ &\left|\left|\xi_{1(\alpha)} - \xi_1\right|\right|^2 + \left|\left|\xi_{2(\alpha)} - \xi_2\right|\right|^2 + \cdots + \left|\left|\xi_{n(\alpha)} - \xi_n\right|\right|^2 < {\epsilon_1}^2, \text{ for all } \alpha \geq \gamma. \\ \Rightarrow\ &\left(\left|\left|\xi_{1(\alpha)} - \xi_1\right|\right|^2 + \left|\left|\xi_{2(\alpha)} - \xi_2\right|\right|^2 + \cdots + \left|\left|\xi_{n(\alpha)} - \xi_n\right|\right|^2\right)^{\frac{1}{2}} < \epsilon_1, \text{ for all } \alpha \geq \gamma. \\ \Rightarrow\ &|||P_\alpha - P|||_R < \epsilon_1, \text{ for all } \alpha \geq \gamma, \text{ where } P = (\xi_1, \xi_2, \cdots, \xi_n). \\ \Rightarrow\ &|||P_\alpha - P|||_R \to 0, \text{ as } \alpha \to \infty. \\ \Rightarrow\ &P_\alpha \to P, \text{ as } \alpha \to \infty, \text{ with respect to real norm.} \end{aligned}$$

Hence, the theorem is proved. □

Theorem 2.27 follows directly from the result established in theorems 2.20 and 2.26.

**Corollary 2.27.** *Every Cauchy sequence in ${\mathbb{C}_2}^n$ converges to an element of ${\mathbb{C}_2}^n$ under complex norm.*

**Theorem 2.28.** *Every Cauchy sequence in ${\mathbb{C}_2}^n$ converges to an element of ${\mathbb{C}_2}^n$ under idempotent norm.*

*Proof.* Let $\{P_k\}_{k=1}^{\infty}$ be an arbitrary Cauchy sequence in ${\mathbb{C}_2}^n$ with respect to idempotent norm, where $P_k = \left(\xi_{1(k)}, \xi_{2(k)}, \cdots, \xi_{n(k)}\right)$. Let $\epsilon > 0$, then there exists $n_0 \in \mathbb{N}$ such that

$$\begin{aligned}
& \left|\left|\left|P_\alpha - P_\beta\right|\right|\right|_{ID} < \sqrt{2}\epsilon, \text{ for all } \alpha, \beta \geq n_0. \\
\Rightarrow\ & \sqrt{2}\left|\left|\left|P_\alpha - P_\beta\right|\right|\right|_{R} < \sqrt{2}\epsilon, \text{ for all } \alpha, \beta \geq n_0, \text{ by theorem 2.23.} \\
\Rightarrow\ & \left|\left|\left|P_\alpha - P_\beta\right|\right|\right|_{R} < \epsilon, \text{ for all } \alpha, \beta \geq n_0.
\end{aligned}$$

Thus, $\{P_k\}_{k=1}^{\infty}$ is the Cauchy sequence in ${\mathbb{C}_2}^n$ with respect to real norm.

Now we consider an $\epsilon_1 > 0$. Then $\dfrac{\epsilon_1}{\sqrt{2}} > 0$ and by theorem 2.26, there exists $P \in {\mathbb{C}_2}^n$ and $n_1 \in \mathbb{N}$ such that

$$\begin{aligned}
& |||P_\alpha - P|||_R < \frac{\epsilon_1}{\sqrt{2}}, \text{ for all } \alpha \geq n_1. \\
\Rightarrow\ & \sqrt{2}\,|||P_\alpha - P|||_R < \epsilon_1, \text{ for all } \alpha \geq n_1. \\
\Rightarrow\ & |||P_\alpha - P|||_{ID} < \epsilon_1, \text{ for all } \alpha \geq n_1, \text{ by theorem 2.23.} \\
\Rightarrow\ & |||P_\alpha - P|||_{ID} \to 0, \text{ as } \alpha \to \infty. \\
\Rightarrow\ & P_\alpha \to P, \text{ as } \alpha \to \infty, \text{ with respect to idempotent norm.}
\end{aligned}$$

Hence, the theorem is proved. □

Theorem 2.29 follows directly from the result established in theorems 1.8 and 2.25 to 2.27.

**Corollary 2.29.** *The structures* $({\mathbb{C}_2}^n, +, \cdot, \times, |||\bullet|||_R)$ *and* $({\mathbb{C}_2}^n, +, \cdot, \times, |||\bullet|||_C)$ *form a modified Banach algebra, and the structures* $({\mathbb{C}_2}^n, +, \cdot, \langle\langle \bullet, \bullet\rangle\rangle_R)$ *and* $({\mathbb{C}_2}^n, +, \cdot, \langle\langle \bullet, \bullet\rangle\rangle_C)$ *are Hilbert spaces.*

Theorem 2.30 follows directly from the result established in theorems 1.8, 2.25 and 2.28.

**Corollary 2.30.** *The structures* $({\mathbb{C}_2}^n, +, \cdot, \times, |||\bullet|||_{ID})$ *and* $({\mathbb{C}_2}^n, +, \cdot, \langle\langle \bullet, \bullet\rangle\rangle_{ID})$ *form a Banach algebra and a Hilbert space, respectively.*

**Definition 2.31.** Let $A = [\xi_{ij}]_{m\times n}$, $B = [\eta_{ij}]_{m\times n}$, and $C = [\zeta_{ij}]_{m\times n}$ be the arbitrary bicomplex matrices of order $m \times n$ such that

$$\begin{aligned}
A &= {}^R\!A_1 + i_1\, {}^R\!A_2 + i_2\, {}^R\!A_3 + i_1 i_2\, {}^R\!A_4 = {}^C\!A_1 + i_2\, {}^C\!A_2 = A^- e_1 + A^+ e_2, \\
B &= {}^R\!B_1 + i_1\, {}^R\!B_2 + i_2\, {}^R\!B_3 + i_1 i_2\, {}^R\!B_4 = {}^C\!B_1 + i_2\, {}^C\!B_2 = B^- e_1 + B^+ e_2, \\
C &= {}^R C_1 + i_1\, {}^R C_2 + i_2\, {}^R C_3 + i_1 i_2\, {}^R C_4 = {}^C C_1 + i_2\, {}^C C_2 = C^- e_1 + C^+ e_2, \\
\xi_{ij} &= {}^{ij}x_1 + i_1\, {}^{ij}x_2 + i_2\, {}^{ij}x_3 + i_1 i_2\, {}^{ij}x_4 = {}^{ij}z_1 + i_2\, {}^{ij}z_2 = (\xi_{ij})^-\, e_1 + (\xi_{ij})^+\, e_2, \\
\eta_{ij} &= {}^{ij}y_1 + i_1\, {}^{ij}y_2 + i_2\, {}^{ij}y_3 + i_1 i_2\, {}^{ij}y_4 = {}^{ij}w_1 + i_2\, {}^{ij}w_2 = (\eta_{ij})^-\, e_1 + (\eta_{ij})^+\, e_2, \\
\zeta_{ij} &= {}^{ij}a_1 + i_1\, {}^{ij}a_2 + i_2\, {}^{ij}a_3 + i_1 i_2\, {}^{ij}a_4 = {}^{ij}u_1 + i_2\, {}^{ij}u_2 = (\zeta_{ij})^-\, e_1 + (\zeta_{ij})^+\, e_2.
\end{aligned}$$

Where $i = 1, 2, \cdots, m$ and $j = 1, 2, \cdots, n$. Then, we define three maps $f_7, f_8$ and $f_9$ as follows:

$$f_7, f_8, f_9 : {\mathbb{C}_2}^{m\times n} \times {\mathbb{C}_2}^{m\times n} \longrightarrow \mathbb{C}_1 \text{ such that}$$

$$f_7(A, B) = \sum_{i=1}^{m}\sum_{j=1}^{n}\sum_{k=1}^{4} ({}^{ij}x_k)\overline{({}^{ij}y_k)}, \tag{22}$$

$$f_8(A, B) = \sum_{i=1}^{m}\sum_{j=1}^{n}\sum_{k=1}^{2} ({}^{ij}z_k)\overline{({}^{ij}w_k)}, \tag{23}$$

$$f_9(A, B) = \sum_{i=1}^{m}\sum_{j=1}^{n} (\xi_{ij})^-\ \overline{(\eta_{ij})^-} + (\xi_{ij})^+\ \overline{(\eta_{ij})^+}. \tag{24}$$

**Theorem 2.32.** *The maps* $f_7, f_8,$ *and* $f_9$*, of the types* (22), (23), *and* (24) *satisfy the following conditions: For all* $l \in \{7, 8, 9\}$,

(a) $f_l(A, A) \geq 0$, *for all* $A \in {\mathbb{C}_2}^{m\times n}$.
(b) $f_l(A, A) = 0$ *if and only if* $A = 0$.

*(c)* $\overline{f_l(B,A)} = f_l(A,B)$*, for all* $A,B \in {\mathbb{C}_2}^{m\times n}$.
*(d)* $f_l(A+B,C) = f_l(A,C)+f_l(B,C)$*, for all* $A,B,C \in {\mathbb{C}_2}^{m\times n}$.
*(e)* $f_l(\alpha A,B) = \alpha f_l(A,B)$*, for all* $A,B \in {\mathbb{C}_2}^{m\times n}$ *and* $\alpha \in \mathbb{C}_1$.

*Proof.* Since

$$\begin{aligned} f_7(A,B) &= \sum_{i=1}^{m}\sum_{j=1}^{n}\sum_{k=1}^{4}({}^{ij}x_k)\overline{({}^{ij}y_k)} \\ &= \sum_{i=1}^{m}\sum_{j=1}^{n}({}^{ij}x_1)\overline{({}^{ij}y_1)} + ({}^{ij}x_2)\overline{({}^{ij}y_2)} + ({}^{ij}x_3)\overline{({}^{ij}y_3)} + ({}^{ij}x_4)\overline{({}^{ij}y_4)} \\ &= Tr\left({}^{R}A_1\ [\overline{({}^{R}B_1)}]^T\right) + Tr\left({}^{R}A_2\ [\overline{({}^{R}B_2)}]^T\right) + Tr\left({}^{R}A_3\ [\overline{({}^{R}B_3)}]^T\right) + Tr\left({}^{R}A_4\ [\overline{({}^{R}B_4)}]^T\right). \end{aligned}$$

In the same manner,

$$\begin{aligned} f_8(A,B) &= Tr\left({}^{C}A_1\ [\overline{({}^{C}B_1)}]^T\right) + Tr\left({}^{C}A_2\ [\overline{({}^{C}B_2)}]^T\right), \\ f_9(A,B) &= Tr\left(A^-\ [\overline{B^-}]^T\right) + Tr\left(A^+\ [\overline{B^+}]^T\right). \end{aligned}$$

The proof now follows from these results. The remaining part is routine and is therefore omitted. Hence, the theorem is proved. □

**Definition 2.33 (Real, Complex, and Idempotent Inner Products on** ${\mathbb{C}_2}^{m\times n}$**:).** Let $f_7, f_8$, and $f_9$ be the maps of the types (22), (23), and (24). Then, by theorems 1.8 and 2.32, $f_7, f_8$, and $f_9$ define inner products on ${\mathbb{C}_2}^{m\times n}$. We refer to these inner products induced by $f_7, f_8$, and $f_9$ as the real inner product, complex inner product, and idempotent inner product on ${\mathbb{C}_2}^{m\times n}$, respectively and denote them by $\langle\bullet,\bullet\rangle_R$, $\langle\bullet,\bullet\rangle_C$, and $\langle\bullet,\bullet\rangle_{ID}$, respectively; that is,

$$f_7(A,B) = \langle A,B\rangle_R, f_8(A,B) = \langle A,B\rangle_C, \text{ and } f_9(A,B) = \langle A,B\rangle_{ID}.$$

The norms induced by the real, complex, and idempotent inner products are called real, complex, and idempotent norms on ${\mathbb{C}_2}^{m\times n}$, respectively. They are denoted by $||\bullet||_R$, $||\bullet||_C$, and $||\bullet||_{ID}$, respectively; that is,

$$||A||_R = \sqrt{\langle A,A\rangle_R}, ||A||_C = \sqrt{\langle A,A\rangle_C}, \text{ and } ||A||_{ID} = \sqrt{\langle A,A\rangle_{ID}}.$$

**Remark 2.34.** It is evident that

$$\begin{aligned} ||A||_R &= \left\{\sum_{i=1}^{m}\sum_{j=1}^{n}\left(||\xi_{ij}||_R\right)^2\right\}^{\frac{1}{2}} = \left\{\sum_{i=1}^{m}\sum_{j=1}^{n}\left(||\xi_{ij}||_C\right)^2\right\}^{\frac{1}{2}} = \left\{\sum_{i=1}^{m}\sum_{j=1}^{n}\left(||\xi_{ij}||\right)^2\right\}^{\frac{1}{2}}, \\ &\quad \text{by theorems 2.4 and 2.8.} \\ &= ||A||_C. \\ ||A||_{ID} &= \left\{\sum_{i=1}^{m}\sum_{j=1}^{n}\left(||\xi_{ij}||_{ID}\right)^2\right\}^{\frac{1}{2}} \\ &= \sqrt{2}\left\{\sum_{i=1}^{m}\sum_{j=1}^{n}\left(||\xi_{ij}||\right)^2\right\}^{\frac{1}{2}}, \text{ by theorem 2.11.} \end{aligned}$$

Thus, $||A||_{ID} = \sqrt{2}\left\{\sum_{i=1}^{m}\sum_{j=1}^{n}\left(||\xi_{ij}||\right)^2\right\}^{\frac{1}{2}} = \sqrt{2}||A||_R = \sqrt{2}||A||_C$.

**Definition 2.35 (Real, Complex, and Idempotent Inner Products on** $\mathbb{C}_2[\xi]_n$**:).** Let $\mathbb{C}_2[\xi]_n$ represent the set of all bicomplex polynomials of degree at most $n$; that is,

$$\begin{aligned} \mathbb{C}_2[\xi]_n &= \{P(\xi) : P(\xi) \text{ is a bicomplex polynomial of degree at most } n\} \\ &= \{P(\xi) : P(\xi) = a_0 + a_1\xi + \cdots + a_n\xi^n, \text{ such that } a_0, a_1, \cdots, a_n \in \mathbb{C}_2\} \end{aligned}$$

Let $P(\xi)$, $Q(\xi)$, and $R(\xi)$ be the arbitrary elements of $\mathbb{C}_2[\xi]_n$ such that

$$\begin{aligned}
P(\xi) &= a_0 + a_1\xi + \cdots + a_n\xi^n, \\
a_k &= {}^kx_1 + i_1\ {}^kx_2 + i_2\ {}^kx_3 + i_1i_2\ {}^kx_4 = {}^kz_1 + i_2\ {}^kz_2 = (a_k)^-\ e_1 + (a_k)^+\ e_2, \\
Q(\xi) &= b_0 + b_1\xi + \cdots + b_n\xi^n, \\
b_k &= {}^ky_1 + i_1\ {}^ky_2 + i_2\ {}^ky_3 + i_1i_2\ {}^ky_4 = {}^kw_1 + i_2\ {}^kw_2 = (b_k)^-\ e_1 + (b_k)^+\ e_2, \\
R(\xi) &= c_0 + c_1\xi + \cdots + c_n\xi^n, \\
c_k &= {}^ka_1 + i_1\ {}^ka_2 + i_2\ {}^ka_3 + i_1i_2\ {}^ka_4 = {}^ku_1 + i_2\ {}^ku_2 = (c_k)^-\ e_1 + (c_k)^+\ e_2,
\end{aligned}$$

where $k = 0, 1, \cdots, n$. We define the algebraic operations of addition $(+)$ and scalar multiplication $(\cdot)$ on $\mathbb{C}_2[\xi]_n$ as follows:

$$\begin{aligned}
P(\xi) + Q(\xi) &= (a_0 + b_0) + (a_1 + b_1)\xi + \cdots + (a_n + b_n)\xi^n, \\
\alpha \cdot P(\xi) &= (\alpha a_0) + (\alpha a_1)\xi + \cdots + (\alpha a_n)\xi^n.
\end{aligned}$$

Then, it is a routine matter to verify that the structure $(\mathbb{C}_2[\xi]_n, +, \cdot)$ forms a vector space over $\mathbb{C}_1$. We define three maps $f_{10}, f_{11}$ and $f_{12}$ as follows:

$$\begin{aligned}
&f_{10}, f_{11}, f_{12} : \mathbb{C}_2[\xi]_n \times \mathbb{C}_2[\xi]_n \longrightarrow \mathbb{C}_1 \text{ such that} \\
&f_{10}(P(\xi), Q(\xi)) = \sum_{k=0}^{n} {}^kx_1\ {}^ky_1 + {}^kx_2\ {}^ky_2 + {}^kx_3\ {}^ky_3 + {}^kx_4\ {}^ky_4, \\
&f_{11}(P(\xi), Q(\xi)) = \sum_{k=0}^{n} {}^kz_1\ \overline{({}^kw_1)} + {}^kz_2\ \overline{({}^kw_2)}, \\
&f_{12}(P(\xi), Q(\xi)) = \sum_{k=0}^{n} (a_k)^-\ \overline{(b_k)^-} + (a_k)^+\ \overline{(b_k)^+}.
\end{aligned}$$

Thus, $f_{10}, f_{11}$, and $f_{12}$ define inner products on $\mathbb{C}_2[\xi]_n$. We refer to these inner products induced by $f_{10}, f_{11}$, and $f_{12}$ as the real inner product, complex inner product, and idempotent inner product on $\mathbb{C}_2[\xi]_n$, respectively and denote them by $\langle \bullet, \bullet \rangle_R$, $\langle \bullet, \bullet \rangle_C$, and $\langle \bullet, \bullet \rangle_{ID}$, respectively; that is,

$$\begin{aligned}
f_{10}(P(\xi), Q(\xi)) &= \langle P(\xi), Q(\xi) \rangle_R, \\
f_{11}(P(\xi), Q(\xi)) &= \langle P(\xi), Q(\xi) \rangle_C, \\
f_{12}(P(\xi), Q(\xi)) &= \langle P(\xi), Q(\xi) \rangle_{ID}.
\end{aligned}$$

## Acknowledgments

The authors sincerely thank their institutions and colleagues for their valuable discussions and constant support throughout this research work.

## Conflict of Interest

The authors declare that there is no conflict of interest regarding the publication of this paper.

*Email address*: prabhatphilosopher@gmail.com

Department of Applied Mathematics, Gautam Buddha University, Greater Noida-201312, India

*Email address*: fahedzul@gmail.com

Department of Mathematics, School of sciences NIT Andhra pradesh,Tadepalligudem, India

*Email address*: amitujlayan@gbu.ac.in

Department of Applied Mathematics, Gautam Buddha University, Greater Noida-201312, India

*Email address*: anjalisharma773@gmail.com

Department of Mathematics, Institute of Applied Sciences, Mangalayatan University, Aligarh-202146, India

*Email address*: akhil.sharma140@gmail.com

Department of Mathematics, Aligarh Muslim University, Aligarh, Uttar Pradesh 202002,India

*Email address*: gyounusg@gmail.com

Department of Mathematical Sciences Islamic University of Science and Technology Kashmir